\documentclass[11pt]{article}
\usepackage[reqno]{amsmath}
\usepackage{amssymb, theorem, enumerate}

\theoremstyle{change}
{\theorembodyfont{\slshape}
\newtheorem{theorem}{Theorem.}[section]
\newtheorem{lemma}[theorem]{Lemma.}
\newtheorem{corollary}[theorem]{Corollary.}
\newtheorem{proposition}[theorem]{Proposition.}
\theorembodyfont{\rmfamily}

\def\proof{\noindent{\textsl{Proof. }}}
\def\sqr#1#2{{\vbox{\hrule height.#2pt
    \hbox{\vrule width.#2pt height#1pt \kern#1pt
        \vrule width.#2pt}\hrule height.#2pt}}}
\def\eqed{\sqr53}
\def\qed{%
    \ifmmode\eqno\eqed
    \else\nobreak\ \hfill\eqed\medbreak\fi}

\def\ints{{\mathbb Z}}
\def\reals{{\mathbb R}}
\def\cx{{\mathbb C}}
\def\cA{{\mathcal A}}

\def\cE{{\mathcal E}}
\def\cH{{\mathcal H}}
\def\cP{{\mathcal P}}
\def\spn{\mathop{\mathrm{span}}\nolimits}
\def\th{\theta}
\def\a{\alpha}
\def\b{\beta}
\def\g{\gamma}
\def\G{\Gamma}
\def\e{\epsilon}
\def\mod#1{{\ \left(\mathrm{mod\ #1}\right)}}

\def\hp{\widehat{p}}
\def\hX{\widehat{X}}
\def\hA{\widehat{A}}
\def\fp{\widetilde{p}}
\def\fa{\widetilde{a}}
\def\fb{\widetilde{b}}
\def\fc{\widetilde{c}}
\def\fX{\widetilde{X}}
\def\fA{\widetilde{A}}
\def\fG{\widetilde{\Gamma}}
\def\fhp{{\mathcal P}}
\def\fhX{{\mathcal X}}
\def\fhA{{\mathcal A}}

\def\one{{\bf 1}}

\newcommand{\ium}{{instantaneous uniform mixing }}
\newcommand{\BM}{{Bose-Mesner }}

\title{Complex Hadamard Matrices, Instantaneous Uniform Mixing and Cubes}
\author{Ada Chan\\
Department of Mathematics and Statistics\\
York University, Toronto, Ontario, Canada M3J 1P3\\
ssachan@yorku.ca\\
}
\begin{document}
\maketitle

\begin{abstract}
We study the continuous-time quantum walks on graphs in the adjacency algebra of
the $n$-cube and its related distance regular graphs.    

For $k\geq 2$, we find graphs in  the adjacency algebra of $(2^{k+2}-8)$-cube
that admit instantaneous uniform mixing at time $\pi/2^k$ and graphs that have perfect state transfer at time $\pi/2^k$.

We characterize the folded $n$-cubes, the halved $n$-cubes and the folded halved $n$-cubes
whose adjacency algebra contains a complex Hadamard matrix.   We obtain the same conditions
for the characterization of these graphs admitting instantaneous uniform mixing.
\end{abstract}

\section{Introduction}

The continuous-time quantum walk on a graph $X$ is given by the transition operator
\begin{equation*}
e^{-itA} = \sum_{k\geq 0} \frac{(-it)^k}{k!}A^k,
\end{equation*}
where $A$ is the adjacency matrix of $X$.
For example, if $X$ is the complete graph on two vertices, $K_2$, then 
\begin{eqnarray*}
e^{-itA} &=& \left(1-\frac{t^2}{2!}+\frac{t^4}{4!}-\ldots \right) I - i\left(t-\frac{t^3}{3!} + \frac{t^5}{5!}-\ldots\right) A\\
&=& \begin{pmatrix} \cos t & -i\sin t \\-i \sin t & \cos t\end{pmatrix}.
\end{eqnarray*}

Being the quantum analogue of the random walks on graphs,  there is a lot of research interest on quantum walks
for the development of quantum algorithms.   Moreover, quantum walks are proved to be universal for quantum computations
\cite{MR2507892}.  In this paper, we focus on the continuous-time quantum walks introduced by Farhi and Gutmann in \cite{MR1638221}.
Please see \cite{MR2852516} and \cite{quant-ph0303081} for surveys on quantum walks.

Since $A$ is real and symmetric,  the operator $e^{-itA}$ is unitary.
We say the continuous-time quantum walk on $X$ is {\sl \ium at time $\tau$} if
\begin{equation*}
|e^{-i\tau A}|_{a,b} = \frac{1}{\sqrt{|V(X)|}},
\qquad \text{for all vertices $a$ and $b$.}
\end{equation*}
This condition is equivalent to  $\sqrt{|V(X)|} e^{-i\tau A}$ being a complex
Hadamard matrix.   Thus if $X$ admits \ium then its adjacency algebra
contains a complex Hadamard matrix.
In $K_2$, the continuous-time quantum walk is \ium at time $\pi/4$.

In \cite{MR2047028},  Moore and Russell discovered that the continuous-time quantum walk on the $n$-cube
is \ium at time $\pi/4$ which is faster than its classical analogue.    
Ahmadi et al. \cite{MR2023606} showed that the complete graph $K_q$ admits \ium if and only if $q\in \{2,3,4\}$.
Best et al. \cite{arXiv:0808.2382} proved that
\ium occurs in graphs $X$ and $Y$ at time $\tau$ if and only if \ium occurs in 
their Cartesian product at the same time.
They concluded that the Hamming graph $H(n,q)$, which is the Cartesian product of $n$ copies of $K_q$,
has \ium if and only if $q\in \{2,3,4\}$.
In the same paper, they also proved that a folded $n$-cube admits \ium if and only if $n$ is odd.

In this paper, we give a necessary condition for the Bose-Mesner algebra of a symmetric association scheme to
contain a complex Hadamard matrix.   
Applying this  condition, we generalize the result of Best et al. to show that the adjacency algebra of $H(n,q)$
contains the adjacency matrix of a graph that admits \ium if and only if $q\in \{2,3,4\}$.
We characterize the halved $n$-cubes and the folded halved $n$-cubes that have \ium.
We obtain the same characterization for  the folded $n$-cubes, the halved $n$-cubes and the folded halved $n$-cubes to have a complex
Hadamard matrix in their adjacency algebras.

A {\sl cubelike graph} is a Cayley graph of the elementary abelian group $\ints_2^d$.
The graphs appear in this paper are distance regular cubelike graphs.
For $k\geq 2$, we find graphs in the adjacency algebra of $H(2^{k+2}-8,2)$ that admit
\ium at time $\pi/2^k$.
Hence, for all $\tau >0$, there exists graphs  that admit \ium at time less than $\tau$.


Given a graph $X$, we use $A(X)$ to denote its adjacency matrix,
and $X_r$ to denote the graph on the vertex set $V(X)$ in which two vertices are adjacent if
they are at distance $r$ in $X$.
We use $I_v$ and $J_v$ to denote the $v\times v$ identity matrix and the $v\times v$
matrix of all ones, respectively.
We drop the subscript if the order of the matrices is clear.

\section{A Necessary Condition}
\label{Section_NecCond}

The graphs we study in this paper are distance regular.   
The adjacency algebra of a distance regular graph is the \BM algebra of a symmetric association scheme.
In this section, we give a necessary condition for a \BM algebra to contain a complex
Hadamard matrix.   This condition is also necessary for a \BM algebra to contain the adjacency matrix
of a graph that admits instantaneous uniform mixing.

A {\sl symmetric association scheme} of order $v$ with $d$ classes is a set 
\begin{equation*}
\cA =\{A_0, A_1, \ldots, A_d\}
\end{equation*}
of $v\times v$ symmetric $01$-matrices satisfying
\begin{enumerate}
\item
$A_0=I$.
\item
$\sum_{j=0}^d A_j = J$.
\item
$A_j A_k = A_kA_j$,  for $j,k=0,\ldots,d$.
\item
$A_j A_k  \in \spn \cA$,  for $j,k=0,\ldots,d$.
\end{enumerate}
For example, if $X$ is a distance regular graph with diameter $d$, then
the set
$\{I, A(X_1), A(X_2),\ldots, A(X_d)\}$ is a symmetric association scheme.

The {\sl \BM algebra} of an association scheme $\cA$ is the span of $\cA$ over $\cx$.
It is known \cite{MR1002568} that the Bose-Mesner algebra contains another basis $\{E_0, E_1, \ldots, E_d\}$
 satisfying
\begin{enumerate}[(a)]
\item
\label{Eqn_Ej1}
$E_j E_k = \delta_{j,k}E_j$,  for $j,k=0,\ldots,d$, and
\item
\label{Eqn_Ej2}
$\sum_{j=0}^d E_j = I$.
\end{enumerate}
Now there exist complex numbers $p_r(s)$'s such that
\begin{equation}
\label{Eqn_Ej4}
A_r = \sum_{s=0}^d p_r(s) E_s,
\qquad \text{for}\ r=0,\ldots,d.
\end{equation}
It follows from Condition~(\ref{Eqn_Ej1}) that
\begin{equation*}
A_r E_s = p_r(s) E_s,
\qquad \text{for $r,s =0,\ldots,d$}.
\end{equation*}
We call the $p_r(s)$'s {\sl the eigenvalues of the association schemes}.
Since the matrices in $\cA$ are symmetric,  the $p_r(s)$'s  are real.

A $v\times v$ matrix $W$ is {\sl type II} if, for $a,b=1,\ldots,v$,
\begin{equation}
\label{Eqn_TypeII}
\sum_{c=1}^v  \frac{W_{ac}}{W_{bc}} = 
\begin{cases}
v & \text{if $a=b$,}\\
0 & \text{otherwise.}
\end{cases}
\end{equation}
A {\sl complex Hadamard matrix} is a type II matrix whose entries have absolute value one.

\begin{proposition}\label{Prop_TypeII}
Let $\cA=\{A_0,A_1,\ldots,A_d\}$ be a symmetric association scheme.
Let $t_0,\ldots,t_d \in \cx\backslash\{0\}$.
The matrix $W=\sum_{j=0}^d t_j A_j$ is type II if and only if
\begin{equation*}
\left[\sum_{h=0}^d p_h(s) t_h \right]\left[\sum_{j=0}^d p_j(s) \frac{1}{t_j}\right] = v,
\qquad \qquad \text{for $s=0,1,\ldots,d$}.
\end{equation*}
\end{proposition}
\proof
The matrix $W$ is type II if and only if
\begin{equation*}
\left[ \sum_{h=0}^d t_h A_h \right]
\left[ \sum_{j=0}^d \frac{1}{t_j} A_j \right] = vI.
\end{equation*}
It follows from Equation~(\ref{Eqn_Ej4}) and Condition~(\ref{Eqn_Ej2}) that
\begin{equation*}
\left[\sum_{h=0}^d \sum_{l=0}^d t_h p_h(l) E_l \right]
\left[ \sum_{j=0}^d \sum_{k=0}^d \frac{1}{t_j} p_j(k) E_k \right] 
= v \sum_{r=0}^d E_r.
\end{equation*}
By Condition~(\ref{Eqn_Ej1}), the left-hand side becomes
\begin{equation*}
\sum_{r=0}^d [\sum_{h=0}^d t_h p_h(r)] [\sum_{j=0}^d \frac{1}{t_j} p_j(r) ] E_r,
\end{equation*}
multiplying $E_s$ to both sides yields  the equations of this proposition.
\qed

Finding type II matrices in the \BM algebra of a symmetric association scheme amounts to solving
the system of equations in Proposition~\ref{Prop_TypeII}, which is not easy as $d$ gets large.
When we limit the scope of the search to complex Hadamard matrices, we get the following
necessary condition which can be checked efficiently.

\begin{proposition}\label{Prop_NecComHad}
If the \BM algebra of $\cA$ contains a complex Hadamard matrix, then
\begin{equation*}
v \leq \left[\sum_{r=0}^d |p_r(s)| \right]^2,
\qquad \text{for $s=0,1,\ldots,d$}.
\end{equation*}
\end{proposition}
\proof
Suppose $W = \sum_{j=0}^d t_j A_j$ is a complex Hadamard matrix.
By Proposition~\ref{Prop_TypeII}, for $s=0,\ldots,d$,
\begin{equation*}
v = \sum_{r=0}^d p_r(s)^2 + \sum_{0\leq h < j \leq d} (\frac{t_h}{t_j}+\frac{t_j}{t_h})p_h(s)p_j(s).
\end{equation*}
Since $|\frac{t_h}{t_j}|=1$, we have $| \frac{t_h}{t_j}+\frac{t_j}{t_h}|\leq 2$ and 
\begin{equation*}
v \leq \sum_{r=0}^d |p_r(s)|^2 + \sum_{0\leq h < j \leq d} 2 |p_h(s)p_j(s)|
= \left[\sum_{r=0}^d |p_r(s)| \right]^2.
\end{equation*}
\qed

Suppose $A(X)$ belongs to the \BM algebra of $\cA$.
If \ium occurs in $X$ at time $\tau$ then $\sqrt{v} e^{-i\tau A(x)}$
is a complex Hadamard matrix and the eigenvalues of $\cA$ satisfy the inequalities in 
Proposition~\ref{Prop_NecComHad}.
For example, the association scheme $\{I_q, J_q-I_q\}$ has eigenvalues $p_0(1)=1$ and $p_1(1)=-1$.
By Proposition~\ref{Prop_NecComHad}, if the
adjacency algebra of $K_q$ contains a complex Hadamard matrix then $q \leq 4$.
Hence \ium does not occur in $K_q$, for $q\geq 5$.

\begin{proposition}\label{Prop_IUMeqn}
Let $X$ be a graph whose adjacency matrix belongs to the \BM algebra of $\cA$.
Let $\th_0,\ldots,\th_d$ be the eigenvalues of $A(X)$ satisfying
\begin{equation*}
A(X)  = \sum_{s=0}^d \th_s E_s.
\end{equation*}
The continuous-time quantum walk of $X$ is \ium at time $\tau$ if and only if
there exist scalars $t_0,\ldots,t_d$ such that 
\begin{equation*}
|t_0|=\ldots =|t_d|=1
\end{equation*}
and
\begin{equation*}
\sqrt{v}  e^{-i\tau \th_s} = \sum_{j=0}^d p_j(s)t_j,
\qquad \text{for}\ s=0,\ldots,d.
\end{equation*}
\end{proposition}
\proof
It follows from Condition~(\ref{Eqn_Ej1}) that
$A(X)^k = \sum_{s=0}^d \theta_s^k E_s$, for $k\geq 0$.
Therefore,
\begin{equation}
\label{Eqn_PropIUMeqn}
\sqrt{v} e^{-i\tau A(X)} = \sqrt{v} \sum_{s=0}^d e^{-i\tau \th_s}E_s
\end{equation}
belongs to $\spn \cA$,
and there exists $t_0, \ldots, t_d$ such that
\begin{equation*}
\sqrt{v} e^{-i\tau A(X)} = \sum_{j=0}^d t_j A_j.
\end{equation*}
By Equation~(\ref{Eqn_Ej4}), we get
\begin{equation*}
\sqrt{v}  e^{-i\tau \th_s} = \sum_{j=0}^d p_j(s)t_j,
\qquad \text{for}\ s=0,\ldots,d.
\end{equation*}
Lastly, $\sqrt{v} e^{-i\tau A(X)}$ is a complex Hadamard matrix
exactly when 
\begin{equation*}
|t_0|=\ldots=|t_d|=1.
\end{equation*}
\qed

For $n,q \geq 2$, 
the {\sl Hamming graph} $H(n,q)$ is the Cartesian product of $n$ copies of $K_q$.
Equivalently, the  vertex set $V$ of the Hamming graph $H(n,q)$ is the set of words of length $n$
over an alphabet of size $q$, and two words are adjacent
if they differ in exactly one coordinate.
The Hamming graph is a distance regular graph on $q^n$ vertices with diameter $n$.
For $j=1,\ldots, n$, $X_j$ is the graph with vertex set $V$ where two
vertices are adjacent when they differ in exactly $j$ coordinates.
Let $A_0=I$ and $A_j=A(X_j)$, for $j=1,\ldots,n$.
Then $\cH(n,q)= \{A_0, A_1,\ldots, A_n\}$ is a symmetric association scheme, 
called the {\sl Hamming scheme}.
For more information on Hamming scheme, please see \cite{MR1002568} and \cite{arXiv:1011.1044}.

It follows from Equation~(4.1) of \cite{arXiv:1011.1044} that
\begin{equation*}
\sum_{j=0}^n x^jA_j = [I_q+x(J_q-I_q)]^{\otimes n},
\end{equation*}
and the eigenvalues of $\cH(n,q)$ satisfy
\begin{equation}
\label{Eqn_KrawSum}
\sum_{j=0}^n p_j(s) x^j = \left(1+(q-1)x\right)^{n-s}(1-x)^s,
\quad \text{for $s=0,\ldots,n$.}
\end{equation}
Using $[x^k] g(x)$ to denote the coefficient of $x^k$ in a polynomial $g(x)$,
we have
for $r,s = 0,\ldots,n$,
\begin{eqnarray}
\nonumber
p_r(s) &=& [x^r] \left(1+(q-1)x\right)^{n-s}(1-x)^s\\
&=& \sum_{h}(-q)^h(q-1)^{r-h}\binom{n-h}{r-h}\binom{s}{h}.
\label{Eqn_Kraw2}
\end{eqnarray}
We get the following characterization from \cite{MR2047028} and \cite{arXiv:0808.2382}.
\begin{theorem}
\label{Thm_Hamming}
The Hamming graph $H(n,q)$ admits \ium if and only if $q\in \{2,3,4\}$.
\end{theorem}

We see from Proposition~\ref{Prop_IUMeqn}  that 
whether a graph $X$ admits \ium depends on only the spectrum of $X$ and the eigenvalues of the \BM algebra containing $A(X)$.
A {\sl Doob graph} $D(m_1,m_2)$ is a Cartesian product of $m_1$ copies
of the Shrihkande graph and $m_2$ copies of $K_4$.   
It is a distance regular graph with the same parameters as the Hamming graph $H(2m_1+m_2,4)$, see Section~9.2B of \cite{MR1002568}.
Since \ium occurs in $H(n,4)$ for all $n\geq 1$, we see that the Doob graph $D(m_1,m_2)$ admits
\ium for all $m_1, m_2 \geq 1$.

\begin{corollary}
The \BM algebra of $\cH(n,q)$ contains a complex Hadamard matrix if and only if $q\in \{2,3,4\}$.
\end{corollary}
\proof
It follows from Equation~(\ref{Eqn_KrawSum}) that
\begin{equation*}
p_r(n)= (-1)^r\binom{n}{r}.
\end{equation*}
By Proposition~\ref{Prop_NecComHad}, if the \BM algebra of $\cH(n,q)$ contains a complex Hadamard matrix, then 
\begin{equation*}
q^n \leq \left[\sum_{r=0}^n |p_r(n)| \right]^2= 4^n.
\end{equation*}
Hence $q \in \{2,3,4\}$.

The converse follows directly from Theorem~\ref{Thm_Hamming}.
\qed

We conclude that 
if $A(X)$ belongs to the \BM algebra of $\cH(n,q)$, for $q\geq 5$, then
\ium does not occur in $X$.

\section{The Cubes}

The Hamming graph $H(n,2)$ is also called the {\sl $n$-cube}.
It is a distance regular graph on $2^n$ vertices
with intersection numbers  
\begin{equation*}
a_j=0,  \quad b_j = (n-j) \quad \text{and}\quad c_j=j, \qquad\text{for $j=0,\ldots,n$}.
\end{equation*}   
It is both bipartite and antipodal, see Section~9.2 of \cite{MR1002568} for details.

It follows from Equation~(\ref{Eqn_KrawSum}) that the eigenvalues of $\cH(n,2)$ satisfy
\begin{equation}
\label{Eqn_pij}
p_r(n-s) = (-1)^rp_r(s) 
\qquad \text{and} \qquad
p_{n-r}(s) = (-1)^sp_r(s),
\end{equation}
for $r,s=0,\ldots,n$.

The proof of Lemma~\ref{Lem_Cube2} uses the following equations, which are Propositions~2.1(3) and 2.3 of \cite{MR1028893}.
\begin{proposition}\label{Prop_CS}
The eigenvalues of $\cH(n,2)$ satisfy
\begin{enumerate}[(a)]
\item
\label{Eqn_CS1}
$p_r(s+1)-p_r(s)=-p_{r-1}(s+1)-p_{r-1}(s)$
and
\item
\label{Eqn_CS2}
$p_{r-1}(s)-p_{r-1}(s+2) = 4\sum_h (-2)^h\binom{n-2-h}{r-2-h}\binom{s}{h}$.
\qed
\end{enumerate}

\end{proposition}

Note that the Kronecker product of two complex Hadamard matrices is a complex Hadamard matrix.
Hence for $\e\in\{-1,1\}$,  
\begin{equation*}
\left[ I_2 + \e i(J_2-I_2) \right]^{\otimes n}
=\sum_{j=0}^n (\e i)^jA_j 
\end{equation*}
is a complex Hadamard matrix in the \BM algebra of $\cH(n,2)$.

Suppose $A(X)$ belongs to the \BM algebra of $\cH(n,2)$ and 
\begin{equation*}
A(X) E_s = \th_s E_s, \quad \text{for $s=0,\ldots,n$}.
\end{equation*}
It follows from Equations~(\ref{Eqn_PropIUMeqn}) and (\ref{Eqn_KrawSum}) that
\begin{equation*}
\sqrt{2^n} e^{-i\tau A(X)} = e^{i\b} \left[ I_2 + \e i(J_2-I_2) \right]^{\otimes n}
\end{equation*}
if and only if
\begin{eqnarray*}
\sqrt{2^n}e^{-i\tau \th_s} &=& e^{i\b} (1+\e i)^{n-s}(1-\e i)^s \\
&=& \sqrt{2^n}e^{i\b}e^{\e i\pi(n-2s)/4},
\qquad\text{for $s=0,\ldots,n$}.
\end{eqnarray*}
This system of equations holds exactly when
\begin{equation*}
e^{i\b} = e^{-i\tau\th_0-\e i \pi n/4}
\end{equation*}
and
\begin{equation*}
e^{-i\tau(\th_s-\th_0)} = e^{-\e i \pi s/2},
\qquad
\text{for}\ s=0,\ldots,n.
\end{equation*}

\begin{lemma}\label{Lem_Cube1}
Suppose $A(X)$ belongs to the \BM algebra of $\cH(n,2)$
and  $A(X) E_s = \th_s E_s$, for $s=0,\ldots,n$.
If there exists $\e \in \{-1,1\}$ such that
\begin{equation*}
\th_s -\th_0 \equiv \e s 2^{k-1} \mod{2^{k+1}},
\qquad \text{for}\ s=0,\ldots,n,
\end{equation*}
then 
there exists $\b \in \reals$ such that
\begin{equation*}
\sqrt{2^n}e^{-i\frac{\pi}{2^k} A(X)} = e^{i\b} [I_2+\e i (J_2-I_2)]^{\otimes n}.
\end{equation*}
That is,
$X$ admits \ium at time $\pi/2^{k}$.
\qed
\end{lemma}

\begin{lemma}\label{Lem_Cube2}
Let $r\geq 1$.
Let $\a$ be the largest integer such that $\binom{n-1}{r-1}$ is divisible by $2^{\a}$.
Suppose $2^{\a+1-h}$ divides $\binom{n-2-h}{r-2-h}$, for $h=0,\ldots \a$.
Then $X_r$ admits \ium at time $\pi/2^{\a+2}$.

Further, if $n$ is even and $r$ is odd then \ium also occurs in $X_{n-r}$ at time $\pi/2^{\a+2}$.
\end{lemma}
\proof
Since $2^{\a+3}$ divides the right-hand side of Proposition~\ref{Prop_CS}~(\ref{Eqn_CS2}), we have
\begin{equation*}
p_{r-1}(s+2)\equiv p_{r-1}(s) \mod{2^{\a+3}},
\qquad \text{for}\ s=0,\ldots, n-2.
\end{equation*}
Applying this congruence repeatedly gives, for $s=0,\ldots,n-1$,
\begin{equation*}
-p_{r-1}(s+1)-p_{r-1}(s)\equiv -p_{r-1}(1)-p_{r-1}(0) \mod{2^{\a+3}}.
\end{equation*}
Now $-p_{r-1}(1)-p_{r-1}(0) = -2\binom{n-1}{r-1}$, which is divisible by $2^{\a+1}$ but not by $2^{\a+2}$.   
Hence there exists $\e = \pm 1$ such that
\begin{equation*}
-p_{r-1}(1)-p_{r-1}(0)\equiv \e 2^{\a+1} \mod{2^{\a+3}}, 
\end{equation*}
and
\begin{equation*}
-p_{r-1}(s+1)-p_{r-1}(s)\equiv \e 2^{\a+1} \mod{2^{\a+3}},
\end{equation*}
for $s=0,\ldots,n-1$.

By Proposition~\ref{Prop_CS}~(\ref{Eqn_CS1}), we have
\begin{equation*}
p_r(s+1)-p_r(s) \equiv \e 2^{\a+1} \mod{2^{\a+3}},
\end{equation*}
and therefore
\begin{equation}
\label{Eqn_Lem_Cube2A}
p_r(s)-p_r(0) \equiv \e s 2^{\a+1} \mod{2^{\a+3}},
\qquad \text{for}\ s=0,\ldots, n.
\end{equation}
By Lemma~\ref{Lem_Cube1}, $X_r$ admits \ium at time $\pi/2^{\a+2}$.

Suppose $n$ is even and $r$ is odd.
By Lemma~\ref{Lem_Cube1}, it is sufficient to show 
\begin{equation*}
p_{n-r}(s)-p_{n-r}(0) \equiv (-1)^{\frac{n+2}{2}}\e s2^{\a+1} \mod{2^{\a+3}},
\qquad \text{for}\ s=0,\ldots, n.
\end{equation*}

When $s$ is even, $2^{\a+2}$ divides $s2^{\a+1}$ and 
$(-1)^{(n+2)/2}\e s 2^{\a+1} \equiv \e s 2^{\a+1} \mod{2^{\a+3}}$.
Applying Equations~(\ref{Eqn_pij}) and (\ref{Eqn_Lem_Cube2A}), we have
\begin{eqnarray*}
p_{n-r}(s)-p_{n-r}(0) &=& p_r(s)-p_r(0)\\
& \equiv &(-1)^{\frac{n+2}{2}}\e s 2^{\a+1} \mod{2^{\a+3}}.
\end{eqnarray*}

When $s$ is odd, Equation~(\ref{Eqn_pij}) gives
$p_{n-r}(s)-p_{n-r}(0) = -p_r(s)-p_r(0)$.
Applying Equation~(\ref{Eqn_Lem_Cube2A}) with $s=1$, we get
\begin{eqnarray}
\nonumber
p_r(1)-p_r(0) &=& \frac{-2r}{n}\binom{n}{r} \\
&\equiv& \e 2^{\a+1} \mod{2^{\a+3}}.
\label{Eqn_Lem_Cube2B}
\end{eqnarray}
Since $r$ is odd, $2^{\a+1}$ divides $\frac{2}{n}\binom{n}{r}$.

If $n\equiv 0 \mod{4}$, then $2^{\a+3}$ divides $2\binom{n}{r}=2p_r(0)$ and
\begin{eqnarray*}
p_{n-r}(s)-p_{n-r}(0) &=& -[p_r(s)-p_r(0)]-2p_r(0)\\
&\equiv& -[p_r(s)-p_r(0)] \mod{2^{\a+3}}\\
&\equiv& (-1)^{\frac{n+2}{2}}\e s 2^{\a+1} \mod{2^{\a+3}}.
\end{eqnarray*}

Suppose $n \equiv 2 \mod{4}$.   By Equation~(\ref{Eqn_Kraw2}), 
\begin{equation*}
2p_r(s) = \sum_j (-1)^j 2^{j+1} \binom{n-j}{r-j}\binom{s}{j}.
\end{equation*}
The hypothesis of this lemma ensures that $2^{\a+3}$ divides $2^{j+1} \binom{n-j}{r-j}\binom{s}{j}$ for $j\geq 2$.
Thus
\begin{equation*}
2p_r(s) \equiv 2\binom{n}{r}-2^2\binom{n-1}{r-1}s \mod{2^{\a+3}}.
\end{equation*}
We see from Equation~(\ref{Eqn_Lem_Cube2B}) that
$2^{\a+1}$ is the highest power of $2$ that divides $\frac{2r}{n}\binom{n}{r}$.
Since $r$ is odd and $n \equiv 2 \mod{4}$, 
$2^{\a+1}$ is the largest power of $2$ that divides $\binom{n}{r}$.
Using our assumption on $\binom{n-1}{r-1}$, 
\begin{equation*}
2p_r(s) \equiv 2^{\a+2}(\g_1-\g_2) \mod{2^{\a+3}},
\end{equation*}
for some odd integers $\g_1$ and $\g_2$.
Therefore, $2p_r(s)$ is divisible by $2^{\a+3}$ and 
\begin{eqnarray*}
p_{n-r}(s) -p_{n-r}(0) 
&=& [p_r(s) - p_r(0)] - 2p_r(s)\\
& \equiv& (-1)^{\frac{n+2}{2}}\e s 2^{\a+1}  \mod{2^{\a+3}}.
\end{eqnarray*}
It follows from Lemma~\ref{Lem_Cube1} that \ium occurs in $X_{n-r}$ at time $2^{\a+2}$.
\qed

To find the $n$'s and $r$'s that satisfy the condition in Lemma~\ref{Lem_Cube2},
we need the following results from number theory,  due to Lucas and Kummer, respectively  (see Chapter~IX of \cite{MR0245499}).
\begin{theorem}
\label{Thm_Lucas}
Let $p$ be a prime.
Suppose the representation of $N$ and $M$ in base $p$ are $n_k\ldots n_1 n_0$ and $m_k\ldots m_1m_0$,
respectively.

Then
\begin{equation*}
\binom{N}{M} \equiv  \binom{n_k}{m_k}\ldots\binom{n_0}{m_0} \mod{p}.
\end{equation*} 
\qed
\end{theorem}
\begin{theorem}
\label{Thm_Kummer}
Let $p$ be a prime.
The largest integer $k$  such that $p^k$  divides $\binom{N}{M}$ is 
the number of carries in the addition of $N-M$ and $M$ in base $p$ representation.
\qed
\end{theorem}

Let $2^{\a}$ be the highest power of $2$ that divides $\binom{n-1}{r-1}$.
That is, there are exactly $\a$ carries in the addition of $n-r$ and $r-1$ in base $2$ representation.
If both $n$ and $r$ are even, then no carry takes place in the right-most digit.
Therefore, there are exactly $\a$ carries in the addition of $n-r$ and $r-2$ in base $2$ representation.
Similarly, when $n$ is odd and $r$ is even, there
are exactly $\a-1$ carries in the addition of $n-r$ and $r-2$ in base $2$ representation.
In both cases, $2^{\a+1}$ does not divide $\binom{n-2}{r-2}$, 
so the hypothesis of Lemma~\ref{Lem_Cube2} does not hold when $r$ is even.

\begin{corollary}\label{Cor_Cube1}
Suppose $n$ is even.
For each odd $r\geq 3$, if $\binom{n-1}{r-1} \equiv 1 \mod{2}$ then $X_r$ and $X_{n-r}$ admits
\ium at time $\frac{\pi}{4}$.
\end{corollary}
\proof
Now both $n-r$ and $r-2$ are odd, there is at least one carry (in the rightmost digit) in the addition of $n-r$ and $r-2$ in base $2$ representation.
By Theorem~\ref{Thm_Kummer}, $2$ divides $\binom{n-2}{r-2}$.
The result follows from applying Lemma~\ref{Lem_Cube2} with $\a=0$.
\qed

\begin{corollary}\label{Cor_Cube2}
Let $n=2^m(2\beta+1)$ where $\beta \geq 0$ and $m\geq 1$.
For each odd $r$ satisfying $1\leq r < 2^m$,
$X_r$ and $X_{n-r}$ admit \ium at time $\pi/4$.
\end{corollary}
\proof
Let $r$ be an odd integer between $1$ and $2^m$.
In base $2$ representation, let $(n-1)$ and $(r-1)$ be $v_k\ldots v_0$ and $u_k\ldots u_0$, respectively.
Then 
$v_j=1$ for $j\leq m-1$ and $u_h=0$ for $h\geq m$, so $\binom{v_j}{u_j}=1$ for all $j$.
By Lucas' Theorem, we have
\begin{equation*}
\binom{n-1}{r-1} 
\equiv 1 \mod{2}.
\end{equation*}
The result follows from Corollary~\ref{Cor_Cube1} and Theorem~\ref{Thm_Hamming}.
\qed

We are now ready to 
show the existence of graphs that admit \ium earlier than time $\pi/4$.
\begin{theorem}
\label{Thm_Cube}
For $k \geq 2$,  $X_{2^{k+1}-7}$, $X_{2^{k+1}-5}$, $X_{2^{k+1}-3}$ and $X_{2^{k+1}-1}$ in $\cH(2^{k+2}-8,2)$ admit 
\ium at time $\pi/2^{k}$.
\end{theorem}
\proof
Let $n=2^{k+2}-8$ and $r=\frac{n}{2}-1$.
Then 
\begin{equation*}
n-r =  2^{k+1}-3 = 2^{k}+2^{k-1}+\ldots+1\cdot 2^3+1\cdot 2^2+0\cdot 2^1+1\cdot 2^0
\end{equation*}
and
\begin{equation*}
r-1 =  2^{k+1}-6 = 2^{k}+2^{k-1}+\ldots+1\cdot 2^3+0\cdot 2^2+1\cdot 2^1+0\cdot 2^0.
\end{equation*}
There are $(k-2)$ carries in the addition of $n-r$ and $r-1$ in base $2$ representation.
By Kummer's Theorem, the highest power of $2$ that divides $\binom{n-1}{r-1}$ is $2^{k-2}$.

We want to show that $2^{k-1-h}$ divides $\binom{n-2-h}{r-2-h}$, for $0 \leq h\leq k-2$.
When $h=0$,
\begin{equation*}
r-2 =  2^{k+1}-7 = 2^{k}+2^{k-1}+\ldots + 1\cdot 2^3+0\cdot 2^2+0\cdot 2^1+1\cdot 2^0,
\end{equation*}
so there are $(k-1)$ carries in the addition of $n-r$ and $r-2$ in base $2$ representation.
By Kummer's Theorem, $2^{k-1}$ divides $\binom{n-2}{r-2}$.

Similarly, there are $(k-2)$ carries in the addition of $n-r$ and $r-3$ in base $2$ representation, so
$2^{k-2}$ divides $\binom{n-3}{r-3}$.


As $h$ increments by 1, the number of $1$'s in the leftmost $(k-2)$ digits in the 
base $2$ representation of $r-2-h$ decreases by at most one.   
Hence there are at least $k-1-h$ carries in the addition of $n-r$ and $r-2-h$ in base $2$ representation,
and $2^{k-1-h}$ divides $\binom{n-2-h}{r-2-h}$, for $h=0,\ldots,k-2$.

It follows from Lemma~\ref{Lem_Cube2} with $\a=k-2$
that $X_{2^{k+1}-5}$ and $X_{2^{k+1}-3}$ admit \ium at time $\pi/{2^k}$.

A similar analysis shows that \ium occurs in $X_{2^{k+1}-1}$ and $X_{2^{k+1}-7}$ at the same time.
\qed

\section{Perfect State Transfer}

Let $u$ and $w$ be distinct vertices in $X$. 
We say that {\sl perfect state transfer} occurs from $u$ to $w$ in the continuous-time quantum walk on $X$ at time $\tau$ if
\begin{equation*}
|(e^{-i\tau A(X)})_{u,w}|=1.
\end{equation*}
We say that $X$ is {\sl periodic} at $u$ with period $\tau$ if
\begin{equation*}
|(e^{-i\tau A(X)})_{u,u}|=1.
\end{equation*}
If $A(X)$ belongs to the \BM algebra of an association scheme and $X$ is periodic at some vertex $u$, then
$X$ is periodic at every vertex.   In this case, we simply say that $X$ is {\sl periodic}.

Consider $X_r$ in the Hamming scheme $\cH(2^m,2)$ when $r$ is odd.
We see from the proof of Corollary~\ref{Cor_Cube2} that
$\binom{2^m-1}{r-1}$ is odd.   It follows from Theorem~2.3 of
\cite{MR2811131} that perfect state transfer occurs in $X_r$
at time $\pi/2$.
Moreover let $1\leq r' \leq 2^m$ be an odd integer distinct from $r$, then
the graph $X_r \cup X_{r'}$ is periodic with period $\pi/2$.

Let $X$ be one of the graphs considered in Corollary~\ref{Cor_Cube2} or
Theorem~\ref{Thm_Cube}.    If $X$ admits \ium at time $\tau$ then
\begin{equation*}
e^{-i\tau A(X)} = \frac{e^{i\beta}}{\sqrt{2^n}}
\begin{pmatrix}
1 & \e i\\
\e i & 1
\end{pmatrix}^{\otimes n},
\qquad \text{for some $\beta \in \reals$ and $\e=\pm1$}.
\end{equation*}
Observe that, for $\e, \e' \in \{-1,1\}$,
\begin{equation}
\label{Eqn_PST1}
\begin{pmatrix}
1 & \e i\\
\e i & 1
\end{pmatrix}\begin{pmatrix}
1 & \e'\  i\\
\e'\  i & 1
\end{pmatrix}
=
\begin{cases}
2\begin{pmatrix}
0&\e i\\\e i&0
\end{pmatrix}
& \text{if $\e=\e'$},\\
2\begin{pmatrix}
1&0\\0&1
\end{pmatrix}
& \text{if $\e\not=\e'$}.\\
\end{cases}
\end{equation}
We see that $X$ is periodic at time $2\tau$.

With the help of the following result in number theory, Theorem~1 of \cite{MR1910963}, we find  graphs in $\cH(2^m, 2)$
and $\cH(2^{k+2}-8,2)$ that have perfect state transfer earlier than $\pi/2$.

\begin{theorem}
\label{Thm_CG}
Let $p$ be prime, $n$ and $k$ be positive integers.  If $p^k$ divides $n$ then
\begin{equation*}
\binom{n-1}{s} \equiv (-1)^{s- \lfloor  s/p\rfloor} \binom{n/p-1}{\lfloor s/p \rfloor} \mod{p^k},
\end{equation*}
for $s=0,\ldots,n-1$.
\end{theorem}

\begin{proposition}
\label{Prop_PST1}
For $m\geq 2$,  and for odd integers $r$ and $r'$ satisfying
\begin{equation*}
1\leq r<r' < 2^{m-1}  \qquad \text{or} \qquad 2^{m-1} < r < r' <2^m,
\end{equation*}
perfect state transfer occurs in the graph $X_r \cup X_{r'}$ of $\cH(2^m,2)$ at time $\pi/4$.
\end{proposition}
\proof
Since $A_r$ and $A_{r'}$ commute, we have
\begin{equation*}
e^{-it (A_r+A_{r'})} = e^{-itA_r}e^{-itA_{r'}}.
\end{equation*}
By Corollary~\ref{Cor_Cube2}, Lemma~\ref{Lem_Cube1} and Equation~(\ref{Eqn_PST1}), it suffices to show that
\begin{equation}
\label{Eqn_PropPST1}
p_r(1) - p_r(0) \equiv
\begin{cases}
-2 \mod{8} & \text{for odd integer $r < 2^{m-1}$} \\
2 \mod{8} & \text{for odd integer $r > 2^{m-1}$}.
\end{cases}
\end{equation}

Let $r$ be an odd integer between $2^b$ and $2^{b+1}$ for some $b\leq m-1$.
Let $s_0=r-1$ and $s_i=\lfloor s_{i-1}/2\rfloor$, for $i=1,\ldots, b$.
Applying Theorem~\ref{Thm_CG} repeatedly gives
\begin{equation*}
\binom{n-1}{r-1} \equiv (-1)^{s_0-s_i}\binom{2^{m-i}-1}{s_i} \mod{4},
\qquad  \text{for $1\leq i \leq \min(m-2,b)$}.
\end{equation*}
If $r<2^{m-1}$, we have $b\leq m-2$, $s_b =1$ and
\begin{eqnarray*}
\binom{n-1}{r-1} &\equiv& (-1)^{s_0-s_b}  (2^{m-b}-1) \mod{4}\\
& \equiv& 1 \mod{4}.
\end{eqnarray*}
If $2^{m-1}<r $, we have $b=m-1$, $2\leq s_{m-2} \leq 3$ and 
\begin{eqnarray*}
\binom{n-1}{r-1} &\equiv& (-1)^{s_0-s_{m-2}}\binom{2^2-1}{s_{m-2}} \mod{4}\\
&\equiv& -1 \mod{4}.
\end{eqnarray*}
Equation~(\ref{Eqn_PropPST1}) follows from the fact that $p_r(1)-p_r(0)=-2\binom{n-1}{r-1}$.
\qed

\begin{proposition}
\label{Prop_PST2}
For integer $k\geq 2$,  perfect state transfer occurs in graphs  $X_{2^{k+1}-7}\cup X_{2^{k+1}-5}$
and $X_{2^{k+1}-3}\cup X_{2^{k+1}-1}$ of $\cH(2^{k+2}-8,2)$ at time $\pi/2^k$.
\end{proposition}
\proof
Let $m=(2^{k+2}-8)/8$.
Applying Theorem~\ref{Thm_Cube} and Lemma~\ref{Lem_Cube1} to $\cH(8m,2)$, we have
\begin{eqnarray*}
p_{4m-1}(1)-p_{4m-1}(0) &=& -2\binom{8m-1}{4m-2}\\
& \equiv&  \e  2^{k-1} \mod{2^{k+1}},
\end{eqnarray*}
for some $\e \in\{-1,1\}$.
Consider
\begin{eqnarray*}
p_{4m-3}(1) - p_{4m-3}(0)  &=& -2 \binom{8m-1}{4m-4}\\
&=& \left[1-4\frac{5m}{(4m+3)(2m+1)}\right] (-2) \binom{8m-1}{4m-2}\\
&\equiv&  \left[1-4\frac{5m}{(4m+3)(2m+1)}\right] \e 2^{k-1} \mod{2^{k+1}}.
\end{eqnarray*}
Since $(4m+3)$ and $(2m+1)$ are coprime with $2^{k+1}$,
\begin{equation*}
p_{4m-3}(1) - p_{4m-3}(0)  \equiv \e 2^{k-1} \mod{2^{k+1}}.
\end{equation*}
By Equation~(\ref{Eqn_PST1}), we have
\begin{equation*}
e^{-i\frac{\pi}{2^{k}} (A_{2^{k+1}-7} + A_{2^{k+1}-5})} = e^{i\b}
\begin{pmatrix}
0&\e i\\\e i&0
\end{pmatrix}^{\otimes n},
\quad \text{for some $\b \in \reals$}.
\end{equation*}

Similarly, 
\begin{eqnarray*}
p_{4m+1}(1)-p_{4m+1}(0) &=& -2 \binom{8m-1}{4m}\\
&\equiv& \e' 2^{k-1} \mod{2^{k+1}}
\end{eqnarray*}
for some $\e' \in \{-1,1\}$.
Then 
\begin{eqnarray*}
p_{4m+3}(1)-p_{4m+3}(0)  &=& -2 \binom{8m-1}{4m+2}\\
&=& \left[1-4\frac{3m}{(2m+1)(4m+1)}\right] \left(p_{4m+1}(1)-p_{4m+1}(0)\right)\\
&\equiv&  \e' 2^{k-1} \mod{2^{k+1}},
\end{eqnarray*}
and perfect state transfer occurs in $X_{2^{k+1}-3}\cup X_{2^{k+1}-1}$ at time $\pi/2^k$.

\qed

\section{Halved $n$-Cube}
The $n$-cube $X$ is a connected bipartite graph of diameter $n$.
When $n\geq 2$, $X_2$ has two components, one of which has the set $\cE$ of binary words 
of even weights as its vertex set.
The {\sl halved $n$-cube}, denoted by $\hX$, is the subgraph of $X_2$ induced by $\cE$.
It is a distance regular graph on $2^{n-1}$ vertices with diameter $\lfloor \frac{n}{2}\rfloor$.
The intersection numbers of $\hX$ are
\begin{equation*}
\widehat{a}_j=2j(n-2j), \quad \widehat{b}_j = \frac{(n-2j)(n-2j-1)}{2}
\quad \text{and}\quad
\widehat{c}_j = j(2j-1),
\end{equation*}
for $j=0,\ldots, \lfloor \frac{n}{2}\rfloor$,
and the eigenvalues of $\hX$ are $p_2(0),p_2(1),\ldots, p_2(\lfloor n/2\rfloor)$.

Let $\widehat{\cA}=\{I,\hA_1, \ldots,\hA_{\lfloor n/2\rfloor}\}$
where $\hA_r= A(\hX_r)$.
We use $\hp_r(s)$ to denote the eigenvalues of $\widehat{\cA}$ and let $\hp_{-1}(s)=0$.
Equation~(11) on page~128 of \cite{MR1002568} states that, 
for $r,s = 0,\ldots,\lfloor n/2\rfloor$,
\begin{equation*}
\hp_1(s) \hp_{r}(s) = \widehat{c}_{r+1} \hp_{r+1}(s) + \widehat{a}_r \hp_{r}(s) +\widehat{b}_{r-1}\hp_{r-1}(s).
\end{equation*}
It is straightforward to verify that $\hp_r(s)=p_{2r}(s)$ satisfies these recursions,
so the eigenvalues of $\widehat{\cA}$ are
\begin{equation}
\label{Eqn_Halvedpij}
\hp_r(s) = p_{2r}(s),
\qquad \text{for}\ r,s=0,\ldots,\lfloor \frac{n}{2}\rfloor.
\end{equation}
For more information on the halved $n$-cube, please see Sections~4.2 and 9.2D of \cite{MR1002568}.

When $n=2m+1$, Equation~(\ref{Eqn_KrawSum}) yields
\begin{equation*}
\sum_{h=0}^n p_h(s) i^h = (1+i)^{2m+1-s}(1-i)^s = 2^mi^{m-s}(1+i),
\qquad \text{for}\ s=0,\ldots,n.
\end{equation*}
The real part of this sum is
\begin{eqnarray}
\label{Eqn_RealPart}
\sum_{r=0}^m p_{2r}(s) (-1)^r &=& \sum_{r=0}^m \hp_{r}(s) (-1)^r\\
&=&\begin{cases}
2^m & \text{if $m-s\equiv 0 \mod{4}$ or $m-s\equiv 3 \mod{4}$},\\ 
-2^m & \text{otherwise}.
\end{cases}\nonumber
\end{eqnarray}
By Proposition~\ref{Prop_TypeII},  $\sum_{r=0}^m (-1)^r \hA_r$ is a complex Hadamard matrix.

\begin{theorem}
For $n \geq 3$,  the adjacency algebra of the halved $n$-cube contains a complex Hadamard
matrix if and only if $n$ is odd.
\end{theorem}
\proof
Suppose $n=2m$.
Using Proposition~\ref{Prop_NecComHad}, it is sufficient to show that 
\begin{equation*}
\left[\sum_{r=0}^m |\hp_r(m-1)|\right]^2 < 2^{2m-1},
\qquad \text{for}\ m\geq 2.
\end{equation*}

It follows from Equations~(\ref{Eqn_KrawSum}) and (\ref{Eqn_Halvedpij}) that for $r \geq 0$,
\begin{eqnarray*}
\hp_r(m-1) 
&=& [x^{2r}] (1+x)^{m+1}(1-x)^{m-1} \\
&=&  [x^{2r}](1+2x+x^2)(1-x^2)^{m-1} \\
&=& (-1)^r \left[ \binom{m-1}{r} - \binom{m-1}{r-1}\right].
\end{eqnarray*}
Hence
\begin{equation*}
|\hp_r(m-1)| = 
\begin{cases}
\binom{m-1}{r} -\binom{m-1}{r-1}  & \text{if $0\leq r \leq \frac{m}{2}$}\\
\binom{m-1}{r-1}-\binom{m-1}{r} & \text{if $\frac{m}{2} < r \leq m$}
\end{cases}
\end{equation*}
and
\begin{eqnarray*}
&&\sum_{r=0}^m |\hp_r(m-1)|\\
&=& \sum_{r=0}^{\lfloor\frac{m}{2} \rfloor}\left[\binom{m-1}{r}-\binom{m-1}{r-1}\right]
+\sum_{r=\lfloor\frac{m}{2}\rfloor+1}^m \left[\binom{m-1}{r-1}-\binom{m-1}{r}\right]\\
&=&2\binom{m-1}{\lfloor\frac{m}{2}\rfloor}.
\end{eqnarray*}
A simple mathematical induction on $m$ shows that 
$4\binom{m-1}{\lfloor \frac{m}{2} \rfloor}^2 < 2^{2m-1}$, 
for $m\geq 2$.

When $n$ is odd, $\sum_{r=0}^m (-1)^r\hA_r$ is a complex Hadamard matrix.
\qed

\begin{theorem}\label{Thm_HalvedIUM}
For $n \geq 3$,  the halved $n$-cube admits \ium if and only if $n$ is odd.
\end{theorem}
\proof
From the above theorem,  the halved $n$-cube does not admit \ium 
when $n\geq 4$ is even.

Suppose $n=2m+1$ and $e^{-2i\tau}=\pm i$.
For $s=0,\ldots,m$, we have 
\begin{equation*}
\hp_1(s)=2(m-s)(m-s+1)-m
\end{equation*}
and
\begin{eqnarray*}
e^{-i\tau\hp_1(s)}
&=&
(e^{-2i\tau})^{(m-s)(m-s+1)} e^{i\tau m}\\
&=& 
\begin{cases}
e^{i\tau m} & \text{if $m-s\equiv 0 \mod{4}$ or $m-s\equiv 3 \mod{4}$},\\ 
-e^{i\tau m} & \text{otherwise}.
\end{cases}
\end{eqnarray*}
We see from Equation~(\ref{Eqn_RealPart}) that
\begin{equation*}
2^m e^{-i\tau\hp_1(s)} = e^{i\tau m}\sum_{r=0}^m (-1)^r \hp_r(s),
\qquad \text{for}\ s=0,\ldots,m.
\end{equation*}
Since $|e^{i\tau m}(-1)^r|=1$, it follows from 
Proposition~\ref{Prop_IUMeqn} that $\hX_1$ admits \ium at time $\frac{\pi}{4}$.
\qed
The halved $2$-cube is the complete graph on two vertices and it admits \ium (see \cite{MR2023606}).

When $n\geq 3$, the halved $n$-cube is isomorphic to the cubelike graph of $\ints_2^{n-1}$ with connection set 
\begin{equation*}
C=\left\{{\bf a} : \text{weight of ${\bf a}$ is $1$ or $2$}\right\}.
\end{equation*}
Applying Theorem~2.3 of \cite{MR2811131} to the halved $n$-cube with even $n$, we see that perfect state
transfer occurs from ${\bf a}$ to ${\bf a} \oplus \one$ at time $\pi/2$.    But this graph does not have instantaneous uniform mixing.

\section{Folded $n$-Cube}
Let $\G$ be a distance regular graph on $v$ vertices with diameter $d$ and 
intersection array $\{b_0, b_1, \ldots, b_{d-1}; c_1, \ldots,c_d\}$.
We say $\G$ is {\sl antipodal} if $\G_d$ is a union of complete graph $K_R$'s, for some fixed $R$.
The vertex sets of the $K_R$'s in $\G_d$ form an equitable partition $\cP$ of $\G$
and the quotient graph of $\G$ with respect to $\cP$ is called
the {\sl folded graph} $\fG$ of $\G$.
When $d>2$,
$\fG$ is a distance regular graph on $\frac{v}{R}$ vertices with  diameter $\lfloor \frac{d}{2}\rfloor$,
see Proposition~4.2.2~(ii) of \cite{MR1002568}.
Moreover $\fG$ has intersection numbers $\fa_j=a_j$, $\fb_j=b_j$ and 
$\fc_j=c_j$ for $j=0,\ldots \lfloor \frac{d}{2}\rfloor-1$ and
\begin{equation*}
\fc_{\lfloor \frac{d}{2}\rfloor} = 
\begin{cases}
c_{\lfloor \frac{d}{2}\rfloor} & \text{if $d$ is odd},\\
R c_{\frac{d}{2}} & \text{if $d$ is even.}
\end{cases}
\end{equation*}
From Proposition~4.2.3~(ii) of \cite{MR1002568}, we see that if the eigenvalues of $\G$ are 
$p_1(0)\geq p_1(1)\geq \ldots\geq p_1(d)$,
then $\fG$ has eigenvalues $\fp_1(j)=p_1(2j)$ for $j=0,\ldots,\lfloor \frac{d}{2}\rfloor$.
The eigenvalues for $\fA_j$'s and $A_j$'s satisfy the same recursive relation (Equation~(11) on Page~128 of \cite{MR1002568}) 
for $j=0,\ldots, \lfloor \frac{d}{2}\rfloor$ when $d$ is odd and
for $j=0,\ldots, \frac{d}{2}-1$ when $d$ is even.
When $d$ is even, $\fp_{\frac{d}{2}}(s) = \frac{1}{R}p_{\frac{d}{2}}(2s)$.
Therefore 
\begin{equation}
\label{Eqn_Foldedpij}
\fp_r(s) = 
\begin{cases}
p_r(2s) & \text{if $0\leq r < \lfloor \frac{d}{2}\rfloor$},\\
p_{\lfloor \frac{d}{2}\rfloor}(2s) & \text{if $d$ is odd and $r={\lfloor \frac{d}{2}\rfloor}$},\\
\frac{1}{R}p_{\frac{d}{2}}(2s) & \text{if $d$ is even and $r={\frac{d}{2}}$}.
\end{cases}
\end{equation}

For each vertex ${\bf a}$ in the $n$-cube $X$, $\one \oplus {\bf a}$ is the unique vertex at distance $n$ from ${\bf a}$.
Therefore $X_n$ is a union of $K_2$'s.  
The {\sl folded $n$-cube} $\fX$ has $2^{n-1}$ vertices, diameter $\lfloor \frac{n}{2} \rfloor$, and eigenvalues
\begin{equation}
\label{Eqn_Foldedn}
\fp_r(s) = 
\begin{cases}
[x^r](1+x)^{n-2s}(1-x)^{2s} & \text{if $0\leq r < \lfloor \frac{n}{2}\rfloor$},\\
[x^{\lfloor \frac{n}{2}\rfloor}](1+x)^{n-2s}(1-x)^{2s} & \text{if $n$ is odd and $r={\lfloor \frac{n}{2}\rfloor}$},\\
[x^{\frac{n}{2}}]\frac{1}{2}(1+x)^{n-2s}(1-x)^{2s} & \text{if $n$ is even and $r={\frac{n}{2}}$}.
\end{cases}
\end{equation}

The folded $n$-cube is isomorphic to the graph obtained from an $(n-1)$-cube by adding the perfect matching in which a vertex ${\bf a}$ is
adjacent to $\one \oplus {\bf a}$.
Best et al. proved the following result, see Theorem~1 of \cite{arXiv:0808.2382}.
\begin{theorem}
For $n \geq 3$,  the folded $n$-cube admits \ium if and only if $n$ is odd.
\qed
\end{theorem}
In particular, the adjacency algebra of the folded $n$-cube contains a complex Hadamard matrix when $n$ is odd.
\begin{theorem}
For $n \geq 3$,  the adjacency algebra of the folded $n$-cube contains a complex Hadamard
matrix if and only if $n$ is odd.
\end{theorem}
\proof
Suppose $n=4m$, for some $m\geq 1$.
We have, for $r=0,\ldots, 2m-1$, 
\begin{eqnarray*}
\fp_r(m) 
&=& [x^r](1+x)^{2m}(1-x)^{2m}\\
&=& 
\begin{cases}
(-1)^{\frac{r}{2}}\binom{2m}{\frac{r}{2}} & \text{if $r$ is even,}\\
0 & \text{otherwise}
\end{cases}
\end{eqnarray*}
and
\begin{equation*}
\fp_{2m}(m) =  (-1)^m\frac{1}{2}\binom{2m}{m}.
\end{equation*}
Now
\begin{eqnarray*}
\sum_{r=0}^{2m} |\fp_r(m)|
&=& \sum_{r=0}^{m-1} \binom{2m}{r} + \frac{1}{2}\binom{2m}{m}\\
&=& \frac{1}{2}\left[ \sum_{r=0}^{2m} \binom{2m}{r}\right]\\
&=& 2^{2m-1}.
\end{eqnarray*}
We have $\left[\sum_{s=0}^{2m} |\fp_s(m)| \right]^2 < 2^{4m-1}$.
By Proposition~\ref{Prop_NecComHad}, the adjacency algebra of the folded $4m$-cube does not contain
a complex Hadamard matrix.

Suppose $n=4m+2$.
By Equation~(\ref{Eqn_Foldedn}),
\begin{eqnarray*}
\fp_r(m) 
&=&
\begin{cases}
1 & \text{if $r=0$,}\\
(-1)^{\lfloor \frac{r}{2}\rfloor}2 \binom{2m}{\lfloor \frac{r}{2}\rfloor} & \text{if $1\leq r< 2m$ is odd,}\\
(-1)^{\frac{r}{2}}\left[ \binom{2m}{\frac{r}{2}} - \binom{2m}{\frac{r}{2}-1}\right] & \text{if $2 \leq r\leq 2m$ is even,}\\
(-1)^m\binom{2m}{m} & \text{if $r=2m+1$.}
\end{cases}
\end{eqnarray*}
Now
\begin{eqnarray*}
\sum_{s=0}^{2m+1} |\fp_s(m)|
&=& 1+\sum_{r=0}^{m-1} 2 \binom{2m}{r}+\sum_{r=1}^{m}\left[\binom{2m}{r}-\binom{2m}{r-1}\right]+\binom{2m}{m}\\
&=& 2^{2m}+\binom{2m}{m}.
\end{eqnarray*}
A simple mathematical induction on $m$ shows that $\left[2^{2m}+\binom{2m}{m}\right]^2 < 2^{4m+1}$, for all integer $m\geq 2$.
We conclude that the adjacency algebra of the folded $(4m+2)$-cube does not contain
a complex Hadamard matrix, for $m\geq 2$.

The folded $6$-cube has eigenvalues
\begin{eqnarray*}
p_0(1)=p_0(2)=1, &	 \quad & p_1(1)=-p_1(2)=2,\\
p_2(1)=p_2(2)=-1 & \quad \text{and} \quad &p_3(1)=-p_3(2)=-2.
\end{eqnarray*}
Let $W=\sum_{j=0}^3 t_j\fA_j$ be a type II matrix.
Adding the equations in Proposition~\ref{Prop_TypeII} for $s=1$ and $s=2$ gives
\begin{equation*}
-\left(\frac{t_0}{t_2}+\frac{t_2}{t_0}\right) - 4 \left(\frac{t_1}{t_3}+\frac{t_3}{t_1}\right)=22.
\end{equation*}
The left-hand side is at most ten if  $|t_0|=|t_1|=|t_2|=|t_3|=1$.
Therefore, the adjacency algebra of the folded $6$-cube does not contain a complex Hadamard matrix.
\qed
The folded $2$-cube is the complete graph on two vertices and it admits \ium (see \cite{MR2023606}).

{

\section{Folded Halved $2m$-Cube}
According to Page~141 of \cite{MR1002568}, the halved $2m$-cube $\hX$ is antipodal with antipodal classes of size two
and the folded $2m$-cube $\fX$ is bipartite for $m\geq 2$.
In addition, the folded graph of $\hX$ is isomorphic to the halved graph of $\fX$.
We use $\fhX$ to denoted the folded graph of $\hX$ which is a distance regular graph on
$2^{2m-2}$ vertices with diameter $\lfloor \frac{m}{2} \rfloor$.
Let $\cA_r=A(\fhX_r)$, for $r=0,\ldots,\lfloor \frac{m}{2}\rfloor$. 

By Equations~(\ref{Eqn_Halvedpij}) and (\ref{Eqn_Foldedpij}), the eigenvalues of the folded halved $2m$-cube are
\begin{equation}
\label{Eqn_FHpij}
\fhp_r(s) = 
\begin{cases}
p_{2r}(2s) & \text{if $0\leq r < \lfloor \frac{m}{2}\rfloor$},\\
p_{2\lfloor \frac{m}{2}\rfloor}(2s) & \text{if $m$ is odd and $r={\lfloor \frac{m}{2}\rfloor}$},\\
\frac{1}{2}p_{m}(2s) & \text{if $m$ is even and $r=\frac{m}{2}$}.
\end{cases}
\end{equation}

\begin{theorem}
The adjacency algebra of the folded halved $2m$-cube contains a complex Hadamard
matrix if and only if $m$ is even.
\end{theorem}
\proof
Suppose $m=2u+1$.  Then
\begin{eqnarray*}
\fhp_r(u) 
&=& [x^{2r}](1+2x+x^2)(1-x^2)^{2u}\\
&=&
\begin{cases}
1 & \text{if $r=0$}\\
(-1)^r\binom{2u}{r}+(-1)^{r-1}\binom{2u}{r-1} & \text{if $1\leq r\leq u$.}
\end{cases}
\end{eqnarray*}
Then
\begin{equation*}
\sum_{r=0}^u |\fhp_r(u)| = 1+\sum_{r=1}^u \left[\binom{2u}{r}-\binom{2u}{r-1}\right]=\binom{2u}{u}.
\end{equation*}
Hence
\begin{equation*}
\left[\sum_{r=0}^u |\fhp_r(u)|\right]^2 < \left[\sum_{r=0}^{2u}\binom{2u}{r}\right]^2 = 2^{4u}.
\end{equation*}
By Proposition~\ref{Prop_NecComHad}, the adjacency algebra of the folded halved $(4u+2)$-cube 
does not contain a complex Hadamard matrix.

Suppose $m=2u$.
By Equations~(\ref{Eqn_FHpij}) and (\ref{Eqn_pij}), 
\begin{eqnarray*}
&&\sum_{r=0}^u (-1)^r \fhp_r(s)\\
&=& \sum_{r=0}^{u-1}(-1)^rp_{2r}(2s) + \frac{1}{2}(-1)^u p_{2u}(2s)\\
&=& \frac{1}{2}\sum_{r=0}^{u-1}(-1)^rp_{2r}(2s) + \frac{1}{2}(-1)^u p_{2u}(2s) +
\frac{1}{2}\sum_{r=0}^{u-1}(-1)^r (-1)^{2s} p_{4u-2r}(2s)\\
&=& \frac{1}{2} \sum_{r=0}^{2u}(-1)^rp_{2r}(2s),
\end{eqnarray*}
which is equal to the real part of $\frac{1}{2}\sum_{j=0}^{4u} i^j p_j(2s)$. 
By Equation~(\ref{Eqn_KrawSum}).
\begin{equation}
\label{Eqn_FHeven}
\frac{1}{2}\sum_{j=0}^{4u} i^j p_j(2s)= \frac{1}{2} (1+i)^{4u-2s}(1-i)^{2s}=(-1)^{u-s}2^{2u-1}.
\end{equation}
By Proposition~\ref{Prop_TypeII},  $\sum_{s=0}^u (-1)^s \fhA_s$ is a complex Hadamard matrix.
\qed

\begin{theorem}
The folded halved $2m$-cube admits \ium if and only if $m$ is even.
\end{theorem}
\proof
Suppose $m=2u$ and $e^{-8i\tau}=-1$.
For $s=0,\ldots,u$, 
\begin{equation*}
\fhp_1(s) = 8(u-s)^2-2u
\end{equation*}
and
\begin{equation*}
2^{2u-1}e^{-i\tau\fhp_1(s)}
= 2^{2u-1}(-1)^{(u-s)^2}e^{2iu\tau},
\end{equation*}
which is equal to $e^{2iu\tau}\sum_{r=0}^u (-1)^r \fhp_r(s)$ from Equation~(\ref{Eqn_FHeven}).
By Proposition~\ref{Prop_IUMeqn}, the folded halved $4u$-cube admits \ium at time $\pi/8$.
\qed


\section*{Acknowledgement}
The author would like to thank Chris Godsil,  Natalie Mullin and Aidan Roy for many interesting discussions.

\bibliography{CubesIUM}
\bibliographystyle{acm}
\end{document}